\documentclass{SCAE}
\numberwithin{equation}{section}

\usepackage{bm}

\numberwithin{equation}{section}

\newtheorem{rmk}{Remark}

	\newcommand{\E}{\mathbb{E}}

	\renewcommand{\P}{\mathbb{P}}
	\renewcommand{\a}{\alpha}

	\newcommand{\la}{\lambda}
	
	\newcommand{\de}{\delta}
	
	\newcommand{\si}{\sigma}
	\newcommand{\Ga}{\Gamma}
	\newcommand{\ga}{\gamma}
	
	\newcommand{\Si}{\pmb\Sigma}
	\renewcommand{\o}{\omega}

	\newcommand{\A}{{\bf A}}
	\newcommand{\B}{{\bf B}}
	\newcommand{\D}{{\bf D}}
		\newcommand{\bfE}{{\bf E}}

	\newcommand{\I}{{\bf I}}

	\renewcommand{\L}{{\bf L}}

	\newcommand{\R}{{\bf R}}
	
	\renewcommand{\S}{{\bf S}}
	\newcommand{\T}{{\bf T}}

	\newcommand{\X}{{\bf X}}
	
	\newcommand{\Y}{{\bf Y}}
	\newcommand{\Z}{{\bf Z}}

	\renewcommand{\(}{\left(}
	\renewcommand{\)}{\right)}

	\newcommand{\lb}{\label}

	\newcommand{\toD}{\stackrel{\mathcal{D}}{\to}}
	\newcommand{\toP}{\stackrel{\mathcal{P}}{\to}}

\begin{document}

\Year{2015} %
\Month{January}
\Vol{55} %
\No{1} %
\BeginPage{1} %
\EndPage{XX} %
\AuthorMark{Hu J \& Bai Z. }
\ReceivedDay{Month 00, 2015}
\AcceptedDay{Month 00, 2015}
\DOI{10.1007/s11425-000-0000-0} 

\title{A review of 20 years of naive tests of significance for  high-dimensional mean vectors and covariance matrices}{
In memory of 10 years since the passing of
 Professor Xiru Chen - a great Chinese statistician}

\author[1]{HU Jiang}{}
\author[1]{BAI Zhidong}{Corresponding author}

\address[{\rm1}]{Key Laboratory for Applied Statistics of MOE, Northeast Normal University, Changchun, 130024, P.R.C.}

\Emails{huj156@nenu.edu.cn, baizd@nenu.edu.cn }\maketitle

 {\begin{center}
\parbox{14.5cm}{\begin{abstract}
In this paper, we introduce the so-called naive tests and give a brief review of the new developments. Naive testing methods are easy to understand and perform robustly, especially when the dimension is large. In this paper, we focus mainly on reviewing some naive testing methods for the mean vectors and covariance matrices of high-dimensional populations, and we believe that this naive testing approach can be used widely in many other testing problems. 

\vspace{-3mm}
\end{abstract}}\end{center}}



\renewcommand{\baselinestretch}{1.2}
\begin{center} \renewcommand{\arraystretch}{1.5}
{\begin{tabular}{lp{0.8\textwidth}} \hline \scriptsize
{\bf Citation:}\!\!\!\!&\scriptsize  HU J, BAI Z D. A review of 20 years of naive tests of significance for  high-dimensional mean vectors and covariance matrices. Sci China Math, 2015, 55, doi: 10.1007/s11425-000-0000-0\vspace{1mm}
\\
\hline
\end{tabular}}\end{center}

\baselineskip 11pt\parindent=10.8pt  \wuhao

\section{ Introduction}
Since its proposal by Hotelling (1931) \cite{Hotelling31G}, the Hotelling $T^2$ test has served as a good test used in multivariate analyses for more than eight decades due to its many useful properties: it is uniformly the most powerful of the affine invariant tests for the hypotheses $H_0:\bm \mu=0$ for the one-sample problem and $H_0: \bm\mu_1=\bm\mu_2$ for the two-sample problem. However, it has a fatal defect in that it is not well defined when the dimension is larger than the sample size or the degrees of freedom. As a remedy,  Dempster (1958) \cite{Dempster58H} proposed his non-exact test (NET) to test the hypothesis of the equality of two multivariate population means, that is, the test of locations in the two-sample problem. 
In 1996, Bai and Saranadasa \cite{BaiS96E} further found that Dempster's NET not only serves as a replacement for the Hotelling $T^2$ to test the hypothesis when the number of degrees of freedom is lower than the dimension but is also more powerful than the Hotelling  $T^2$ when the dimension is large, but not too large, such that  $T^2$ is well defined. They also proposed the asymptotic normal test (ANT) to test the same hypothesis and strictly proved that both the NET and ANT have similar asymptotic power functions that are higher than those of the Hoteling $T^2$ test. 
Thus, their work raised an important question that classical multivariate statistical procedures need to re-examine when the dimension is high. To call attention to this problem, they entitled their paper ``The Effect of High Dimension''.

That paper was published nearly 20 years ago and has been cited in other studies more than 100 times to date in Web of Science. It is interesting that more than 95\% of the citations were made in the past 10 years. This pattern reveals that high-dimensional data analysis has attracted much more widespread attention since the year 2005 than it had received previously. In the theory of hypothesis testing, of course, the most preferred test is the uniformly most powerful test. However, such a test does not exist unless the distribution family has the property of a monotone likelihood ratio for which the parameter can only be univariate. Hence, there is no uniformly most powerful test for multivariate analysis. Therefore, the optimal procedure can only be considered for smaller domains of significance tests, such as unbiased tests or invariant tests with respect to specific transformation groups. The Hotelling $T^2$ was derived based on the likelihood ratio principle and proved to be the most powerful invariant test with respect to the affine transformation group (see Page 174 of \cite{Anderson03I}). A serious point, however, is that the likelihood ratio test must be derived under the assumption that the likelihood of the data set exists and is known, except for the unknown parameters. In a real application, it is impossible to verify that the underlying distribution is multivariate normal or has any other known form of the likelihood function.  Thus, we would like to use another approach to set up a test for some given hypothesis: choose $h(\bm\theta)$ as a target function for the hypotheses such that the null hypothesis can be expressed as $h(\bm\theta)=0$ and the alternative as $h(\bm\theta)>0$ and then look for a proper estimator $\bm{\hat\theta}$ of the parameter $\bm\theta$. 
Then, we reject the hypothesis if $h(\bm{\hat\theta})>h_0$ such that $P_{H_0}(h(\bm{\hat\theta})>h_0)=\alpha$.  For example, for the Hotelling test of the difference of two sample means, one can choose $h(\bm \mu_1,\bm\mu_2,\Si)=(\bm\mu_1-\bm\mu_2)'\Si^{-1}(\bm\mu_1- \bm\mu_2)$, the estimators $\hat{\bm\mu}_i=\bar{\X}_i$, $i=1,2$, and $\hat\Si=\S$ for the sample means and sample covariance matrix. Dempster's NET and Bai and Saranadasa's ANT simply use $h(\bm \mu_1,\bm\mu_2)=\|\bm\mu_1-\bm \mu_2\|^2$ and $\hat{\bm\mu}_i=\bar{\X}_i$, $i=1,2$. That is, the Hotelling test uses the squared Mahalanobis distance, whereas the NET and ANT use the squared Euclidean distance. We believe that the reason why the NET and ANT are more powerful for large dimensions than the Hotelling test is because the target function of the latter involves too many nuisance parameters in $\Si$, which cannot be well estimated. Because the new tests focus only on the naive target function instead of the likelihood ratio, we call them the naive tests, especially the ones that are independent of the nuisance parameters, which generally ensures higher power. 

In 1996, Bai and Saranadasa \cite{BaiS96E} raised the interesting point that one might prefer adopting a test of higher power and approximate size rather than a test of exact size but much lower power. The naive tests have undergone rapid development over the past twenty years, especially over the past 10. In this paper, we give a brief review of the newly developed naive tests, which are being applied to a wide array of disciplines, such as genomics, atmospheric sciences, wireless communications, biomedical imaging, and economics. However, due to the limited length of the paper, we cannot review all of the developments and applications in all directions, although some of them are excellent and interesting for the field of high-dimensional data analysis. In this paper, we focus mainly on reviewing some naive testing methods (NTMs) for the mean vectors and covariance matrices of high-dimensional populations.

Based on the NTMs, many test statistics have been proposed for high-dimensional data analysis.   
Throughout this paper, we suppose that there are $k$ populations and that the observations $\X_{i1},\dots,\X_{in_i}$ are  $p$-variate independent and identically distributed (i.i.d.) random sample vectors from the $i$-th population, which have the mean vector $\bm\mu_i$ and the covariance matrix $\Si_i$. Moreover, except where noted, we work with the following model assumptions:
\begin{itemize}
	\item[(A1)] $\X_{ij}:=(X_{ij1},\dots,X_{ijp})'=\bm\Ga_i\Z_{ij}+\bm \mu_i,\quad \mbox{for } i=1,\dots k, j=1\dots,n_i,$ where $\bm \Ga_i$ is a $p\times m$ non-random matrix for some $m\geq p$ such that $\bm\Ga_i\bm\Ga_i'=\Si_i$, and $\{\Z_{ij}\}_{j=1}^{n_i}$ are $m$-variate i.i.d. random vectors satisfying $\E(\Z_{ij})=0$ and $Var(\Z_{ij})=\I_m$, the $m\times m$ identity matrix;
	
	\item[(A2)] $
	\frac{n_i}{n}\to \kappa_i\in(0,1)\quad i=1,\dots k,  \mbox{~as~~} n\to\infty,
	$ where  $n=\sum_{i=1}^kn_i$.

	\end{itemize}
	
	Denote 
\begin{align*}
	\bar{\X}_i=\frac1{n_i}\sum_{j=1}^{n_{i}}\X_{ij}~ \mbox{and}~ \S_i=\frac1{n_i-1}\sum_{j=1}^{n_{i}}(\X_{ij}-\bar{\X}_i)(\X_{ij}-\bar{\X}_i)'=(s^{(i)}_{ij}).
\end{align*}
When $k=1$, the subscripts $i$ or $1$ are suppressed from $n_i$, $n_1$, $\bm\Ga_i$, $\bm\mu_i$ and so on, for brevity. 

Throughout the paper, we denote by $\toP$ the convergence in probability and by $\toD$ the convergence in distribution. 

The remainder of the paper is organized as follows: In Section 2, we review the sample location parameters. In subsection 2.1, we introduce the findings of Bai and Saranadasa \cite{BaiS96E}. In subsection 2.2, we introduce Chen and Qin \cite{ChenQ10T}'s test based on the unbiased estimator of the target function. In subsection 2.3, we review Srivastava and Du's work on the scale invariant NTM, based on the modified component-wise squared Mahalanobis distance. In subsection 2.4, we introduce Cai et al's NTM based on the Kolmogorov distance, i.e., the maximum component of difference. In subsection 2.5, we introduce some works on the extensions to MANOVA and contrast tests, that is, tests for problems of more than two samples.    In Section 3, we introduce some naive tests of hypotheses on covariances. In subsection 3.1, we introduce the naive test proposed by Ledoit  and Wolf \cite{LedoitW02S} on the hypothesis of the one-sample covariance matrix and the spherical test. In subsection 3.2, we introduce the NTM proposed by Li and Chen (2012) \cite{LiC12T}. In subsection 3.3, we introduce Cai's NTM on covariances based on the Kolmogorov distance. We also review the testing of the structure of the covariance matrix in subsection 3.4.  In Section 4, we make some general remarks on the development of NTMs.

\section{Testing the population locations}

\subsection{Asymptotic powers of $T^2$, NET and ANT}

In this section,  we first consider the simpler one-sample problem by NTM. That is, the null hypothesis is  $H_0:\bm \mu_1=\bm \mu_0$.
Under the assumption (A1) with $k=1$, and testing the hypothesis
\begin{align*}
H_0: \bm\mu=\bm\mu_0\quad \mbox{v.s.}\quad  H_1: \bm\mu\neq\bm\mu_0,
\end{align*}  
it is easy to check that $\E \bar{\X}=\bm \mu$. Thus, to set up a test of this hypothesis, we need to choose some norms of the difference $\bm\mu-\bm\mu_0$. There are three types of norms to be chosen in the literature: the Euclidean norm, the Maximum component norm and the Mahalanobis squared norm.  Let us begin from the classical one.

The most famous test  is the so-called Hotelling $T^2$ statistic, 
\begin{align}\lb{reht}
T^2=n(\bar{\X}-\bm\mu_0)'\S^{-1}(\bar{\X}-\bm \mu_0)
\end{align}
which was proposed by Hotelling (1931) \cite{Hotelling31G} and is a natural multi-dimensional extension of the squared univariate Student's t-statistic. If the ${\Z_{j}}$s are normally distributed, the Hotelling $T^2$ statistic is shown to be the likelihood ratio test for this one-sample problem and to have many optimal properties. Details can be found in any textbook on multivariate statistical analysis, such as \cite{Anderson03I,Muirhead82A}.  
It is  easy to verify that $\bar{\X}$ and $\bm \S$ are  unbiased, sufficient and complete estimators of the parameters $\bm\mu$ and $\bm \Si$ and that, as mentioned above, the target function is chosen as the Mahalanobis squared distance of the population mean $\bm\mu$ from the hypothesized mean $\bm \mu_0$, which is also the Euclidean norm of $\Si^{-1/2}(\bm{\mu}-\bm\mu_0)$.  Thus, we can see that the Hotelling $T^2$ statistic is a type of NTM, and we simply need to obtain its (asymptotic) distribution. 
It is well known that under the null hypothesis, $\frac{(n-p)}{p(n-1)}T^2$ has an $F$-distribution with degrees of freedom $p$ and $n-p$, and when  $p$ is fixed, as $n$ tends to infinity, $T^2$ tends to a chi-squared distribution with degrees of freedom $p$.  If we assume $y_{n}=p/n\to y\in(0,1)$ and $\X_{j}$ are normally distributed, following  Bai and Saranadasa \cite{BaiS96E}, we may easily derive that 
\begin{align}\lb{bais1}
\sqrt{\frac{(1-y_{n})^3}{2ny_{n}}}\(T^2-\frac{ny_{n}}{1-y_{n}}
-\frac{n\|\bm\delta\|^2}{1-y_{n}}\)\toD N(0,1), \quad \mbox{as } n\to \infty,
\end{align}
where $\bm \delta=\bm\Si^{-1/2}(\bm\mu-\bm\mu_0)$. By (\ref{bais1}), it is easy to derive that the asymptotic power function of the $T^2$ test satisfies 
$$
\beta_{H}-\Phi\(-\xi_{\alpha}+\sqrt{\frac{n(1-y)}{2y}}\|\bm\delta\|^2\)\to 0.
$$
 Here and throughout the paper, $\Phi$ is used for the distribution function of a standard normal random variable, and $\xi_\a$ is its upper $\a$ quantile. It should be noted that the above asymptotic distribution of the Hotelling $T^2$ statistic \eqref{bais1} still holds without the normality assumption. The details can be found in \cite{PanZ11C}.



Next, we derive the asymptotic power for ANT. In this case, the target function is chosen as $h(\bm\mu)=\|\bm\mu-\bm\mu_0\|^2$, and the natural estimator of $\bm\mu$ is $\bar{\X}=\frac1{n}\sum_{i=1}^{n}\X_{i}$. It is easy to derive that 
\begin{eqnarray}\label{bais11}
\E\|\bar{\X}\|^2&=&\|\bm\mu\|^2+\frac1{n}{\rm tr}\Si\\
Var(\|\bar{\X}\|^2)&=&\frac2{n}{\rm tr}\Si^2+4\bm\mu'\Si\bm\mu\nonumber \\ 
&&+\frac{2}{\sqrt{n}}EZ_{1}^3\sum_{i=1}^m\bm\mu'\bm\gamma_i(\bm\gamma_i'\bm\gamma_i)+\frac1{n}(EZ_{1}^4-3)\sum_{i=1}^m(\bm\gamma_i'\bm\gamma_i)^2\label{bais12}
\end{eqnarray}
where $\bm\gamma_i$ is the $i$-th column of the matrix $\Ga$. Under the conditions
\begin{eqnarray}\label{bais13}
(\bm\mu-\bm\mu_0)'\Si(\bm\mu-\bm\mu_0)&=&o(\frac1{n}{\rm tr}\Si^2),\\
\lambda_{\max}(\Si)&=&o(\sqrt{{\rm tr}\Si^2}),\label{bais14}
\end{eqnarray}
we have 
\[
Var(\|\bar{\X}-\bm\mu_0\|^2)=\left(\frac2{n}{\rm tr}\Si^2+\frac1{n}(EZ_{1}^4-3)\sum_{i=1}^m(\bm\gamma_i'\bm\gamma_i)^2\label{bais12}\right)(1+o(1)).
\]
Under the conditions (\ref{bais13}) and (\ref{bais14}), using the moment method or martingale decomposition method, one can prove that 
\begin{equation}
\frac{\|\bar{\X}-\bm\mu_0\|^2-\E(\|\bar{\X}-\bm\mu_0\|^2)}{\sqrt{Var(\|\bar{\X}-\bm\mu_0\|^2)}}\to N(0,1)
\label{bais15} 
\end{equation} 
To perform the test for the hypothesis $H_0:\bm\mu=\bm\mu_0$ vs. $H_1:\bm\mu\ne \bm\mu_0$, it is necessary to construct ratio-consistent estimators of $\E(\|\bar{\X}-\bm\mu_0\|^2)$ and $Var(\|\bar{\X}-\bm\mu_0\|^2)$ under the null hypothesis. It is obvious that 
$\frac1{n}{\rm tr}\Si$ can be estimated by $\frac1{n}{\rm tr}(\S)$. The variance can be simply estimated by $\frac1{n}\left({\rm tr}(\S^2)-\frac{1}{n}{\rm tr}^2(\S)\right)$ if $EZ_{1}^4=3$. In the general case, it can be estimated by $\frac1{n}\hat\sigma_n^2$, where 
\begin{eqnarray}\label{bais16}
\quad\quad  \quad \hat\sigma_n^2&=&\frac{1}{(n)_5}\sum_{j_1,\cdots,j_5\atop distinct}{\rm tr}\left(({\bf X}_{j_1}-{\bf X}_{j_2})({\bf X}_{j_1}-{\bf X}_{j_3})'({\bf X}_{j_1}-{\bf X}_{j_4})({\bf X}_{j_1}-{\bf X}_{j_5})'\right)\nonumber\\
&&\ \ -\frac{1}{(n)_6}\sum_{j_1,\cdots,j_6\atop distinct}{\rm tr}\left(({\bf X}_{j_1}-{\bf X}_{j_2})({\bf X}_{j_1}-{\bf X}_{j_3})'({\bf X}_{j_6}-{\bf X}_{j_4})({\bf X}_{j_6}-{\bf X}_{j_5})'\right),
\end{eqnarray}
where the summations above are taken for all possibilities that $j_1,\cdots, j_s$, $s=5$ or $6$, distinctly run over $\{1,\cdots,n\}$, and $(n)_l=n(n-1)\cdots(n-l+1)$. 
Using the standard limiting theory approach, one may prove that $\hat\sigma^2_n$ is a ratio-consistent estimator of $\sigma_n^2$, where 
$$
\sigma_n^2=2{\rm tr}\Si^2+(\E Z_{1}^4-3)\sum_{i=1}^m(\bm\gamma_i'\bm\gamma_i)^2\label{bais12}.
$$ 
Therefore, the test rejects $H_0$ if
$$
\|\bar{\X}-\bm\mu_0\|^2>\frac1{n}{\rm tr}(\S)+\frac1{\sqrt{n}}\xi_\alpha \hat\sigma_n.
$$
From this result, it is easy to derive that under conditions (\ref{bais13}) and (\ref{bais14}), the asymptotic power of ANT is
\begin{eqnarray*}
	\beta_{ANT}\simeq\Phi\left(-\xi_\alpha+\frac{\sqrt{n}\|\bm\mu-\bm\mu_0\|^2}{\sqrt{\hat\sigma_n^2}}\right)
	\simeq \Phi\left(-\xi_\alpha+\frac{\sqrt{n}\|\bm\mu-\bm\mu_0\|^2}{\sqrt{\sigma_n^2}}\right).
\end{eqnarray*}
Comparing the expressions of the asymptotic powers of Hotelling test and ANT, one sees that the factor $\sqrt{1-y}$ appears in the asymptotic power of Hotelling's test but not in that of the ANT. This difference shows that the ANT has higher power than does the $T^2$ test when $y$ is close to 1.

Moreover, if $p$, the dimension of the data, is larger than $n-1$, the degrees of freedom, then $T^2$ is not well defined, and there is no way to perform the significance test using it.

\begin{rmk}\label{rmk1}
	In the real calculation of $\hat\sigma_n^2$, the computation using the expression of (\ref{bais16}) is very time consuming. To reduce the computing time, we should rewrite it as
	\begin{eqnarray}\label{bais17}
	\hat\sigma_n^2&=&\frac{1}{n}\sum_{j=1}^{n}({\bf X}_{j}'{\bf X}_{j})^2-\frac{4}{(n)_2}\sum_{j_1,j_2\atop distinct}{\bf X}_{j_1}'{\bf X}_{j_1}{\bf X}_{j_1}'{\bf X}_{j_2}-\frac{1}{(n)_2}\sum_{j_1,j_2\atop distinct}{\bf X}_{j_1}'{\bf X}_{j_1}
	{\bf X}_{j_2}'{\bf X}_{j_2}\nonumber\\
	&&+\frac{6}{(n)_3}\sum_{j_1,j_2,j_3\atop distinct}{\bf X}_{j_1}'{\bf X}_{j_2}{\bf X}_{j_1}'{\bf X}_{j_3}+\frac{2}{(n)_3}\sum_{j_1,j_2,j_3\atop distinct}{\bf X}_{j_1}'{\bf X}_{j_1}{\bf X}_{j_2}'{\bf X}_{j_3}-\frac{4}{(n)_4}\sum_{j_1,j_2,j_3,j_4\atop distinct}{\bf X}_{j_1}'{\bf X}_{j_2}{\bf X}_{j_3}'{\bf X}_{j_4},\quad\quad\quad 
	\end{eqnarray}
	where each summation runs over all possibilities in which the indices involved are distinct. 
	
	It is easy to see that to calculate the estimator $\hat\sigma_n^2$ using (\ref{bais17}) is very time consuming: for example, to calculate the last term, one needs to compute $2pn^4$ multiplications. To further reduce the computation time, one may use the inclusion-exclusion principle to change the last five sums into forms that are easier to calculate. For  example, the last sum $I_6$ can be written as 
	\begin{eqnarray}\label{eqbais10}
		I_6=({\bf X}'{\bf X})^2-2{\bf X}'{\bf X}a-4{\bf X}'{\bf X_{(2)}}{\bf X}+a^2+2{\rm tr}({\bf X}_{(2)}^2)+8{\bf X}'{\bf X}_{(3)}-6b,
	\end{eqnarray}
	where 
	\begin{eqnarray*}
		{\bf X}&=&\sum_{j=1}^{n}{\bf X_{j}}, \ \ \ \ 
		a=\sum_{j=1}^{n}{\bf X}_{j}'{\bf X}_{j}, \ \ \ \
		{\bf X_{(2)}}=\sum_{j=1}^{n}{\bf X}_{j}{\bf X}_{j}',\\
		{\bf X_{(3)}}&=&\sum_{j=1}^{n}{\bf X}_{j}'{\bf X}_{j}{\bf X}_{j}, \ \ \ \ \
		b=\sum_{j=1}^{n}({\bf X}_{j}'{\bf X}_{j})^2.
	\end{eqnarray*}
	Here, the coefficients of various terms can be found by the following arguments: Let $\Omega$ denote the fact that there are no restrictions between the indices $j_1,\cdots,j_4$, and let $A_{ik}$ denote the restriction that $j_i=j_k$, $i<k\le 4$, which is called an equal sign, or an edge between vertices $i$ and $j$. 
	
	The sum $I_6$ in which the indices $j_1,\cdots,j_4$ are distinct can be considered the indices running over the 
set $\prod_{1\le i<k\le 4}(\Omega-A_{ik})$. By expanding the product, one may split the sum $I_6$ into a signed sum of several sums:  the first sum runs over $\Omega$, followed by the subtraction of 6 sums with one equal sign; add 15 sums with two equal signs; subtract 20 sums with three equal signs, and so on, and finally add the sum with all six equal signs. Now, the first one runs over $\Omega$, that is, there are no restrictions among the four vertices $1,2,3,4$, which simply gives the first term in (\ref{eqbais10}). The sum with the equal sign $A_{12}$ is given by 
	$$
	\sum_{A_{12}}{\bf X}_{j_1}'{\bf X}_{j_2}{\bf X}_{j_3}'{\bf X}_{j_4}=\sum_{j=1}^{n}\sum_{j_3=1}^{n}\sum_{j_4=1}^{n_1}{\bf X}_{j}'{\bf X}_{j}{\bf X}_{j_3}'{\bf X}_{j_4}=a{\bf X}'{\bf X};
	$$
Similarly, the sum under the equal sign $A_{34}$ is also $a{\bf X}'{\bf X}$.
These two cases give $-2{\bf X'X}a$ in the second term in (\ref{eqbais10}); the other 4 cases with one equal sign give 
$-4{\bf X}'{\bf X_{(2)}}{\bf X}$ in the third term. For example,
	$$
	\sum_{A_{13}}{\bf X}_{j_1}'{\bf X}_{j_2}{\bf X}_{j_3}'{\bf X}_{j_4}=\sum_{j=1}^{n}\sum_{j_2=1}^{n}\sum_{j_4=1}^{n_1}
	{\bf X}_{j}'{\bf X}_{j_2}{\bf X}_{j}'{\bf X}_{j_4}={\bf X}'{\bf X}_{(2)}{\bf X};
	$$
Under the equal sign $A_{14}$, $A_{23}$ or and $A_{24}$, the sum again has the form ${\bf X}'{\bf X}_{(2)}{\bf X}$.  
By similar arguments, one can show that the sum with the two equal signs $A_{12}$ and $A_{34}$ is given by $a^2$ in the fourth term; the sums with two equal signs $A_{13}$ and $A_{24}$ or $A_{14}$ and $A_{23}$ are given by $2{\rm tr}({\bf X}_{(2)}^2)$ in the fifth term; the sums for the other 12 cases with two equal signs, such as $A_{12}$ and $A_{23}$, are given by $12{\bf X}'{\bf X}_{(3)}$; there are 4 cases in which the three equal signs make three indices equal and leave one index free of the rest (or, equivalently, three edges forming a triangle), which contribute $-4 {\bf X}'{\bf X}_{(3)}$; and two cases give a final contribution of $8{\bf X}'{\bf X}_{(3)}$ in the seventh term  of (\ref{eqbais10}). There are 16 other cases of three equal signs that imply all indices $j_1,\cdots,j_4$ are equal, giving a sum of $b$. Additionally, if there are more than three equal signs, the indices $j_1,\cdots, j_4$ are also all identical; thus, we obtain the coefficient for $b$ as  $-16+{6\choose 4}-{6\choose 5}+1=-6$. Therefore, the splitting (\ref{eqbais10}) is true.  
 
 Similarly, one may show that 
 \begin{eqnarray*}
 I_1&=& b\\
 I_2&=&{\bf X}'{\bf X}_{(3)}-b\\
 I_3&=& a^2-b\\
 I_4&=&{\bf X}'{\bf X}_{(2)}{\bf X}-2{\bf X}'{\bf X}_{(3)}-{\rm tr}({\bf X}_{(2)}^2)+2b\\
 I_5&=&{\bf X}'{\bf X}a-a^2-2{\bf X}'{\bf X}_{(3)}+2b.
 \end{eqnarray*}
 Finally, one can calculate the estimator of the variance by
 \begin{eqnarray*}
 \hat \sigma_n^2&=&\frac{n^2-n+2}{(n-1)_3}b-\frac{4}{(n-2)_2}{\bf X}'{\bf X}_{(3)}-\frac{n^2-3n+4}{(n)_4}a^2+\frac{6n-2}{(n)_4}{\bf X}'{\bf X}_{(2)}{\bf X}-\frac{6n-10}{(n)_4}{\rm tr}({\bf X}_{(2)}^2)\\
 &&+\frac{2n+2}{(n)_4}{\bf X'X}a-\frac{4}{(n)_4}({\bf X'X})^2.
 \end{eqnarray*}
 From the expressions above, the numbers of multiplications to be calculated for the terms above are 
 $n(p+1)$, $np^3$,  $np+1$,   $np^3$, $p^2$, $p+1$ and $p+1$, respectively. Thus, by using this formula, the computation time will be reduced significantly. 
\end{rmk}

Now, consider the two-sample location problem of multivariate normal distributions,
$H_0:\bm\mu_1=\bm\mu_2$, with a common covariance matrix $\Si$. 
The classical test for this hypothesis is the Hotelling $T^2$ test
$$
T^2=\frac{n_1n_2}{n_1+n_2}(\bar{\X}_1-\bar{\X}_2)'{\S}^{-1}(\bar{\X}_1-\bar{\X}_2)
$$
where $\bar {\X}_i=\frac1{n_i}\sum_{j=1}^{n_i}{\X}_{ij}$, $i=1,2$ and 
$${\S}=\frac1{n_1+n_2-2}\left(\sum_{i=1}^2\sum_{j=1}^{n_i}({\X}_{ij}-\bar{\X}_i)({\X}_{ij}-\bar{\X}_i)'\right).
$$
In 1958 and 1960, Dempster published two papers, 
\cite{Dempster58H} and \cite{Dempster60S}, 
in which he argued that if $p$, the dimension of the data, is larger than $N=n_1+n_2-2$, the degrees of freedom, then $T^2$ is not well defined, and there is no way to perform the significance test using $T^2$. Therefore, he proposed the so-called NET (non-exact test) as follows: Arrange the data $\bm{\mathcal X}=(\X_{11},\cdots,\X_{1n_1},\X_{21},\cdots,\X_{2n_2})$ as a $p\times n$ matrix, where $n=n_1+n_2$. Select an $n\times n$ orthogonal matrix $H$ and transform the data matrix to
$\bm{\mathcal Y}=\bm{\mathcal X}H=(\bm y_1,\cdots,\bm y_{n})$ such that 
\begin{eqnarray*}
	\bm y_1&\sim& N\Big(\frac{n_1}{\sqrt{n}}\bm\mu_1+\frac{n_2}{\sqrt{n}}\bm\mu_2,\Si\Big)\\
	\bm y_2&\sim& N\Big(\sqrt{\frac{n_1n_2}{n}}(\bm\mu_1-\bm\mu_2),\Si\Big)\\
	\bm y_3, \cdots,\bm y_{n}&\stackrel{i.i.d.}\sim& N(0,\Si).
\end{eqnarray*}
Then, he defined his NET by 
$$
T_D=\frac{\|\bm y_2\|^2}{\|\bm y_3\|^2+\cdots+\|\bm y_{n}\|^2}.
$$
He claimed that as $n_1$ and $n_2$ increase, using the so-called chi-square approximation of $\|\bm y_j\|^2$ for $i=2,3,\cdots,n$, $$nT_D\sim F_{r,nr},$$ where $N=n_1+n_2-2$ and $r=(\textbf{{\rm tr}}\Sigma)^2/{\rm tr}\Sigma^2$. Comparing $NT_D$ with $T^2$, we find that Dempster simply replaced $\S$ by ${\rm tr} (\S)\I_p$ to smooth the trouble when $\S$ is singular.

Bai and Saranadasa \cite{BaiS96E} observed that Dempster's NET is not only a remedy for the $T^2$ test when it is not well defined but is also more powerful than $T^2$ even when it is well defined, provided that the dimension $p$ is large compared with the degrees of freedom $N$. Based on Dempster's NET, Bai and Saranadasa \cite{BaiS96E} also proposed the so-called ANT (asymptotic normality test) based on the normalization of $\|\bm y_2\|^2$. They established a CLT (central limit theorem) as follows:
$$
\frac{M_n}{\sqrt{Var(M_n)}}\toD N(0,1),
$$ 
where $M_n=\|\bar{\X}_1-\bar{\X}_2\|^2-\frac{n}{n_1n_2}{\rm tr}(\S)$. To perform the significance test for $H_0$, Bai and Saranadasa proposed the ratio-consistent estimator of $Var(M_n)$ by
$$
\widehat{Var(M_n)}=\frac{2(N+2)(N+1)N}{n_1^2n_2^2(N-1)}\({\rm tr}(\S^2)-\frac1{N}{\rm tr}^2{\S}\).
$$

Bai and Saranadasa \cite{BaiS96E} proved that when dimension $p$ is large compared with $n$, both NET and ANT are more powerful than the $T^2$ test, by deriving the asymptotic power functions of the three tests. Under the conditions that $p/n\to y\in (0,1)$ and $n_1/n\to\kappa\in(0,1)$, the power function of the $T^2$ test asymptotically satisfies 
$$
\beta_{H}(\delta)-\Phi\(-\xi_\a+\sqrt{\frac{n(1-y)}{2y}}\kappa(1-\kappa)\|\delta\|^2\)\to 0,
$$  
where $\delta=\Si^{-1/2}(\bm\mu_1-\bm\mu_2)$.

Under the assumption A2 and 
\begin{eqnarray}\label{dem1}
	\bm\mu'\Si\bm\mu&=&o\(\frac{n}{n_1n_2}{\rm tr}(\Si^2)\)\\
	\lambda_{\max}(\Si)&=&o(\sqrt{{\rm tr}\Si^2}),\label{dem2}
\end{eqnarray}
where $\bm\mu=\bm\mu_1-\bm\mu_2$, Bai and Saranadasa \cite{BaiS96E} also proved that the power function of Dempster's NET satisfies
$$
\beta_{D}(\bm\mu)-\Phi\(-\xi_\a+\frac{n\kappa(1-\kappa)\|\bm\mu\|^2}{\sqrt{2{\rm tr}{\Si^2}}}\)\to 0.$$

Without the normality assumption, under the assumptions A1 and A2 and the conditions (\ref{dem1}) and (\ref{dem2}), 
Bai and Saranadasa \cite{BaiS96E} proved that the power function of their ANT has similar asymptotic power to the NET, that is,
$$
\beta_{BS}(\bm\mu)-\Phi\(-\xi_\a+\frac{n\kappa(1-\kappa)\|\bm\mu\|^2}{\sqrt{2{\rm tr}{\Si^2}}}\)\to 0.
$$

\subsection{Chen and Qin's approach} Chen and Qin (2010) \cite{ChenQ10T} argued that the main term of Bai and Saranadasa's ANT contains squared terms of sample vectors that may cause non- robustness of the test statistic against outliers and thus proposed  
an unbiased estimator of the target function $\|\bm\mu_1-\bm\mu_2\|^2$, given by 
$$
T_{CQ}=\frac1{n_1(n_1-1)}\sum_{i\ne j}{\X}_{1i}'{\X_{1j}}
+\frac1{n_2(n_2-1)}\sum_{i\ne j}{\X}_{2i}'{\X_{2j}}-\frac2{n_1n_2}\sum_{i,j}{\X}_{1i}'{\X_{2j}}
$$ 
Chen and Qin  \cite{ChenQ10T} proved that $\E T_{CQ}=\|\bm\mu_1-\bm\mu_2\|^2$, and under the null hypothesis, 
$$
Var(T_{CQ})=\frac{2}{n_1(n_1-1)}{\rm tr}(\Si_1^2)+\frac{2}{n_2(n_2-1)}{\rm tr}(\Si_2^2)+\frac{4}{n_1n_2}{\rm tr}(\Si_1\Si_2)(1+o(1)).
$$
Similarly to Bai and Saranadasa (1996), under the conditions 
\begin{eqnarray}
	\frac{n_1}{n}&\to& \kappa\nonumber\\
	{\rm tr}(\Si_i\Si_j\Si_l\Si_h)&=&o\big({\rm tr}^2(\Si_1+\Si_2)^2\big),  \ \mbox{for } i,j,l,h=1\mbox{ or }2,\label{cq11*}\\
	(\bm\mu_1-\bm\mu_2)'\Si_i(\bm\mu_1-\bm\mu_2)&=&o\big(n^{-1}{\rm tr}(\bm\Si_1+\bm\Si_2)^2\big), \ \mbox{for } i=1\mbox{ or }2,\label{cq111}\\
	\mbox{or}\ \ \ \ n^{-1}{\rm tr}(\bm\Si_1+\bm\Si_2)^2&=&o\big((\bm\mu_1-\bm\mu_2)'\Si_i(\bm\mu_1-\bm\mu_2)\big),\ \mbox{for } i=1\mbox{ or }2,\label{cq112}
\end{eqnarray}
Chen and Qin proved that 
$$
\frac{T_{CQ}-\|\bm\mu_1-\bm\mu_2\|^2}{\sqrt{Var(T_{CQ})}}\stackrel{D}\to N(0,1).
$$
To perform the test for $H_0:\bm\mu_1=\bm\mu_2$ with the target function $h(\bm\mu_1,\bm\mu_2)=\|\bm\mu_1-\bm\mu_2\|^2$, they proposed the estimator for $Var(T_{CQ})$ to be 
$$
\hat\sigma_n^2=\frac{2}{n_1(n_1-1)}\widehat{{\rm tr}(\Si_1^2)}+\frac{2}{n_2(n_2-1)}\widehat{{\rm tr}(\Si_2^2)}+\frac{4}{n_1n_2}\widehat{{\rm tr}(\Si_1\Si_2)},
$$
where
\begin{eqnarray*}
	\widehat{{\rm tr}(\Si_i^2)}&=&\frac1{n_i(n_i-1)}\sum_{j\ne k}{\bf X}_{ik}'({\bf X}_{ij}-\bar{\bf X}_{i(jk)}){\bf X}_{ij}'({\bf X}_{ik}-\bar{\bf X}_{i(jk)}),\\
	\widehat{{\rm tr}(\Si_1\Si_2)}&=&\frac1{n_1n_2}\sum_{j=1}^{n_1}\sum_{k=1}^{n_2}{\bf X}_{2k}'({\bf X}_{1j}-\bar{\bf X}_{1(j)}){\bf X}_{1j}'({\bf X}_{2k}-\bar{\bf X}_{2(k)}),
\end{eqnarray*}
and $\bar{\bf X}_{i(*)}$ denotes the sample mean of the $i$-th sample, excluding the $*$-th vectors, as indicated in the braces.

Applying the central limit theorem, Chen and Qin derived the asymptotic power functions for two cases:  
\begin{eqnarray}
	\beta_{CQ}&\sim& \begin{cases}\Phi\left(-\xi_\alpha+\frac{n\kappa(1-\kappa)\|\bm\mu_1-\bm\mu_2\|^2}{\sqrt{2{\rm tr}(\kappa\Si_1+(1-\kappa)\Si_2)^2}}\right)&\
		\quad \mbox{if (\ref{cq111}) holds,}\cr \Phi\left(\frac{n\kappa(1-\kappa)\|\bm\mu_1-\bm\mu_2\|^2}{\sqrt{2{\rm tr}(\kappa\Si_1+(1-\kappa)\Si_2)^2}}\right)&\
		\quad \mbox{if (\ref{cq112}) holds.}\cr\end{cases}  
\end{eqnarray}
\begin{rmk}
The expression of the asymptotic power under the condition (\ref{cq112}) ((3.5) in Chen and Qin \cite{ChenQ10T}) may contain an error in that the denominator of the quantity inside the function $\Phi$ should be $\sigma_{n2}$ in Chen and Qin's notation, that is, 
$2\sqrt{\frac1{n_1}(\bm\mu_1-\bm\mu_2)'\Si_1(\bm\mu_1-\bm\mu_2)+\frac1{n_2}(\bm\mu_1-\bm\mu_2)'\Si_2(\bm\mu_1-\bm\mu_2)}$, instead of $\sigma_{n1}$. However, the asymptotic power is 1 under the condition (\ref{cq112}). Therefore, the typo does not affect the correctness of the expression of the asymptotic power. This point is shown by the following facts:
  
  By the condition (\ref{cq11*}), we have 
$$
(\bm\mu_1-\bm\mu_2)'\Si_i(\bm\mu_1-\bm\mu_2)\le \lambda_{\max}(\Si_i)\|\bm\mu_1-\bm\mu_2\|^2\le o\Big(\sqrt{{\rm tr}(\Si_i^2)}\|\bm\mu_1-\bm\mu_2\|^2\Big).
$$
Therefore, 
$$
\sigma_{n2}^2\le o\big(\sigma_{n1}\|\bm\mu_1-\bm\mu_2\|^2\big). 
$$
Consequently, 
\begin{eqnarray*}
\frac{\|\bm\mu_1-\bm\mu_2\|^2}{\sigma_{n2}}\ge \frac{\sigma_{n2}}{o(\sigma_{n1})}\ge \sqrt{\frac{\sigma_{n2}^2}{\sigma_{n1}^2}}\ge M
\end{eqnarray*}
for any fixed constant $M$, where the last step follows from the condition (\ref{cq112}). Regarding Chen and Qin's expression, one has 
$$
\frac{\|\bm\mu_1-\bm\mu_2\|^2}{\sigma_{n1}}\ge \frac{\|\bm\mu_1-\bm\mu_2\|^2}{\sigma_{n2}}\ge M.
$$ 

\end{rmk}
For the one-sample location problem, Chen and Qin \cite{ChenQ10T}  modified $T_{BS}$ and proposed  
\begin{align*} 
	T_{CQ}=\frac{1}{n(n-1)}\sum_{i\neq j}X_{i}'X_{j}
\end{align*}
and showed that under the condition ${\rm tr}\Si^4=o({\rm tr}^2\Si^2)$ (which is equivalent to (\ref{bais14})), as $\min\{p,n\}\to\infty$, 
\begin{align}\lb{CQT1}
	\frac{T_{CQ}}{\sqrt{\frac{2}{n(n-1)}{\rm tr}(\sum_{i\neq j}({\X}_{i}-\bar{\X}_{(i,j)})\X_{i}'(\X_j-\bar{\X}_{(i,j)})\X_{j})}}\toD N(0,1).
\end{align}
Note that the difference between the statistics in \eqref{CQT1} and the LHS of \eqref{bais15} with the denominator replaced by the estimator $\frac1{n}({\rm tr}{\S}^2-\frac1{n}{\rm tr}^2{\S})$ is in  the denominators,  that are the estimators of   
${\rm tr}\Si^2$.

\begin{rmk}It has been noted that the main part $T_{CQ}$ of Chen and Qin's test is exactly the same as Bai and Saranadasa's $M_n$ because they are both unbiased estimators of the target function $\|\bm\mu_1-\bm\mu_2\|^2$ and are functions of the complete and sufficient statistics of the mean vectors and covariance matrices for the two samples. We believe that Chen and Qin's idea of an unbiased estimator of the target function helped them propose a better estimator of the asymptotic variance of the test such that their test performed better than did Bai and Saranadasa's ANT in simulation.  In addition, there is an improved statistic for the Chen and Qin  test by thresholding methods, which was recently proposed by Chen et al. \cite{ChenL14T}. Wang et al \cite{WangP15H} also proposed a test for the hypothesis under the elliptically distributed assumption, which can be viewed as a nonparametric extension of $T_{CQ}$. 
\end{rmk}

\subsection{Srivastava and Du's approach}

While acknowledging the defect of the Hotelling test, as indicated in 
\cite{Dempster58H, BaiS96E}, Srivastava and Du (2008) \cite{SrivastavaD08T} noted that the NET and ANT are not scale invariant, which may cause lower power when the scales of different components of the model are very different. Accordingly, they proved the following modification to the ANT:

\begin{align*}
	T_{SD,1}=(\bar{\X}-\bm\mu_0)'\D_{\S}^{-1}(\bar{\X}-\bm\mu_0),
\end{align*}  
where $\D_{\S}=Diag(s_{11},\cdots,s_{pp})$ is
the diagonal matrix of the sample covariance matrix $\S$. Let  ${\mathcal R}$ be  the population correlation matrix.
Then, under the condition that 
\begin{eqnarray*}
	0<\lim_{p\to\infty} \frac{{\rm tr}{\mathcal R^i}}{p}&<&\infty, \ \mbox{ for } i=1,2,3,4;\\
	\lim_{p\to\infty}\frac{\lambda( \mathcal R)}{\sqrt{p}}&=&0,
\end{eqnarray*}
where $\lambda( \mathcal R)$ is the largest eigenvalue of the correlation matrix $\mathcal R$.

Srivastava and Du  \cite{SrivastavaD08T} showed that if $n\asymp p^\eta$ and  $\frac{1}{2}<\eta\leq 1$, 
\begin{align*}
	\frac{nT_{SD,1}-\frac{(n-1)p}{n-3}}{\sqrt{2\({\rm tr}\R^2-\frac{p^2}{n-1}\)c_{p,n}}}\toD N(0,1), \quad \mbox{as } n\to \infty,
\end{align*}
where $\R$ is the sample correlation matrix, i.e., $\R=\D_{\S}^{-1/2}\S\D_{\S}^{-1/2}$ and $c_{p,n}= 1+{\rm tr}\R^2/p^{3/2}$. 
  They showed that 
the asymptotic power of $T_{SD,1}$ under the local alternative, as $n\to \infty$,
\begin{align*}
	\beta_{SD}= P\(\frac{nT_{SD,1}-\frac{(n-1)p}{n-3}}{\sqrt{2\({\rm tr}\R^2-\frac{p^2}{n-1}\)c_{p,n}}}> \xi_{-\a} \)
	\to \Phi\(-\xi_{\a}+\frac{n(\bm\mu-\bm\mu_0)'D_{\Si}^{-1}(\bm\mu-\bm\mu_0)}{\sqrt{2{\rm tr}{\mathcal R}^2}}\),
\end{align*} 
where $\D_{\Si}$ is the diagonal matrix of population covariance matrix $\Si$.
Later, Srivastava \cite{Srivastava09T} modified the asymptotic results above to cases, where the adjusting term $c_{p,n_1}$ in the last test statistic is replaced by 1 and the restriction for $\eta$ is relaxed to $0<\eta\leq 1$. Further, by excluding $\sum_{i=1}^n(\X_{i}-\bm\mu_0)'\D_{\S}^{-1}(\X_{i}-\bm\mu_0)$ from  $T_{SD,1}$ and modifying $\D_{\S}$, Park and Ayyala \cite{ParkA13T} obtained another NTM test statistic:
\begin{align*} 
	T_{PA}=\frac{n-5}{n(n-1)(n-3)}\sum_{i\neq j}\X_{i}'\D_{\S_{(i,j)}}^{-1}\X_{j}
\end{align*}
where $\D_{\S_{(i,j)}}=Diag(s^{((i,j))}_{11},\cdots,s^{((i,j))}_{pp})$ is
the diagonal matrix of the sample covariance matrix excluding the sample points $\X_i$ and $\X_j$, i.e., $S_{(i,j)}=(n-3)^{-1}\sum_{k\neq i,j}(\X_{k}-\bar{\X}_{(i,j)})(\X_{k}-\bar{\X}_{(i,j)})'$ and $\bar{\X}_{(i,j)}=(n-2)^{-1}\sum_{k\neq i,j}\X_{k}$. 

Srivastava and Du also considered the two-sample location problem with the common covariance matrix  $\Si$ \cite{SrivastavaD08T} and proposed the testing statistic
$$
T_{SD,2}=\frac{n_1n_2}{n}(\bar{\X}_1-\bar{\X}_2)'\D_{\S}^{-1}(\bar{\X}_1-\bar{\X}_2)
$$
Under similar conditions for the CLT of $T_{SD,1}$, they proved that
\begin{equation}
	\frac{T_{SD,2}-\frac{Np}{N-2}}{\sqrt{2\left({\rm tr}{\R}^2-\frac{p^2}{n}\right)c_{n,p}}}\toD N(0,1).
	\label{limit1}
\end{equation}
They then further derived the asymptotic power function 
$$
\beta_{SD}(\bm\mu)\sim \Phi\left(-\xi_\a+\frac{\kappa(1-\kappa)\bm\mu'\bm\D_{\Si}^{-1}\bm\mu}{\sqrt{2{\rm tr}{\mathcal R}^2}}\right).
$$
≈\begin{rmk}
	The advantage of this statistic is that the terms $\X_{i},~\D_{S_{(i,j)}}$ and $\X_{j}$ are all independent such that it is easy to obtain the approximation 
	\begin{align*}
		 \E T_{PA}=\bm\mu'\E \D_{\Si}^{-1}\bm\mu \simeq \bm \mu'\D_{\Si}^{-1}\bm \mu,
	\end{align*}
	which is similar to $\E (nT_S-{p(n-1)}/{(n-3)})$, as given in \cite{SrivastavaD08T}. This point shows that both $T_{SD,1}$ and $T_{PA}$ are NTM tests based on the target function $\bm \mu'\D_{\Si}^{-1}\bm \mu$. The idea that they use to exclude the bias ${p(n-1)}/{(n-3)}$ is similar to $T_{CQ}$, which removes the bias estimator ${\rm tr}S$ given in $T_{BS}$. 
	Park and Ayyala also gave the asymptotic distribution of $T_{PA}$ under the null hypothesis, that is,
	\begin{align*}
		\frac{\sqrt{n(n-1)}T_{PA}}{\sqrt{2\widehat{{\rm tr}{\mathcal R}^2}}}\to N(0,1),
	\end{align*} 
	where $\widehat{{\rm tr}{\mathcal R}^2}$ is a ratio-consistent estimator of ${\rm tr}{\mathcal R}^2$, i.e., $\widehat{{\rm tr}{\mathcal R}^2}/{\rm tr}{\mathcal R}^2\toP1$, and  
	\begin{align*}
		\widehat{{\rm tr}{\mathcal R}^2}=\frac{1}{n(n-1)}\sum_{i\neq j}X_{i}'\D_{\S_{(i,j)}}^{-1}(\X_{j}-\bar{\X}_{(i,j)})\X_{j}'\D_{\S_{(i,j)}}^{-1}(\X_{i}-\bar{\X}_{(i,j)}).
	\end{align*} 
	They then showed that the asymptotic power of the test $T_{PA}$ is the same as the asymptotic power of $T_{SD}$.
	Recently, Dong et al. \cite{DongP16S} gave a shrinkage estimator of the diagonals of the population covariance matrix $D_{\Si_1}$ and showed that the shrinkage-based Hotelling test performs better than the unscaled Hotelling test and the regularized Hotelling test when the dimension is large.
\end{rmk}

\begin{rmk}
	For $\Si_1\neq \Si_2 $, Srivastava et al. \cite{SrivastavaK13T} used  $\D=\D_{\S_1}/n_1+\D_{\S_2}/n_2$ instead of $D_{\S}$ in $T_{SD,2} $. 
	For the case in which the population covariance matrices are diagonal, 
	Wu et al. \cite{WuG06M} constructed a statistic by summing up the squared component-wise $t$-statistics for missing data, and Dong et al. \cite{DongP16S} proposed a shrinkage-based diagonalized Hotelling’s test.

	%

\end{rmk}

\subsection{Cai et al's idea}

Cai et al \cite{CaiX14H} noted that all of the NTM tests associated with target functions based on the Euclidean norm or Mahalanobis distance require the condition that 
\begin{eqnarray}
	n\|\bm\mu\|^2/\sqrt{p}\to\infty,
	\label{ntmcond}
\end{eqnarray}
to distinguish the null and alternative hypotheses with probability tending to 1. This condition does not hold if only a few components of $\bm\mu$ have the order $O(1/\sqrt{n})$ and all others are 0. Therefore, they proposed using the $L_\infty$ norm or, equivalently, the Kolmogorov distance.  

Indeed, Cai et al's work compensates for the case where 
\begin{eqnarray}
	\sqrt{n}\max_{i\le p}|\mu_i|\to\infty.
	\label{caicond}
\end{eqnarray}
Note that neither of the conditions (\ref{ntmcond}) and (\ref{caicond}) implies the other. The condition (\ref{caicond}) is weaker than (\ref{ntmcond}) only when the bias vector  $\bm\mu-\bm\mu_0$ or $\bm\mu_1-\bm\mu_2$  is sparse. 

Now, we introduce the work of Cai et al. \cite{CaiL14T}, who developed another NTM test based on the Kolmogorov distance, 
which performs more powerfully against sparse alternatives in high-dimensional settings. 

Supposing that $\Si_1=\Si_2=\Si$ and  $\{\X_1, \X_2\}$ satisfy the sub-Gaussian-type  or polynomial-type tails condition, Cai et al. proposed the test statistic
\begin{align*}
	T_{CLX}=\frac{n_1n_2}{n_1+n_2}\max_{1\leq i\leq p}\left\{\frac{{\mathcal X}_i^2}{\o_{ii}}\right\},
\end{align*}
where $\widehat{\Si^{-1}}(\bar X_1-\bar X_2):=({\mathcal X}_1,\dots,{\mathcal X}_p)'$ and $\widehat{\Si^{-1}}:=\Omega=(\o_{ij})_{p\times p}$ is the constrained
$l_1$-minimization for the inverse matrix estimator of $\Si^{-1}$. Here, the so-called constrained $l_1$-minimization for the inverse matrix estimator is defined by
\begin{eqnarray*}
	\widehat{\Si^{-1}}&=& {\rm arg}\min_{\Omega=(\omega_{ij})}\Big\{ \sum_{ij}|\omega_{ij}|; {\mbox{subject to }} \|{\S}\Omega-\I_p\|_{\infty}\le \ga_n\Big\} 
\end{eqnarray*}
where $\ga_n$ is a tuning parameter, which may generally be chosen as $C\sqrt{{\log p}/n}$ for some large constant $C$. 
For more details on the properties of $l_1$-minimization estimators, the reader is referred to
\cite{CaiL11C}. Under the null hypothesis $H_0$ and some spectrum of population covariance matrix conditions, for any $x\in R$, as $\min\{n,p\}\to \infty$,
\begin{align*}
	\P(T_{CLX}-2\log(p)-\log\log(p)\leq x)\to \mbox{Exp}\(-\frac{1}{\pi}\mbox{Exp}\(-\frac{x}{2}\)\).
\end{align*}
To evaluate the performance of their maximum absolute components test, they also proved the following result:

Suppose that $C^{-1}_0\le \lambda_{\min}(\Si)\le \lambda_{\max}(\Si)\le C_0$ for some constant $C_0>1$; $k_p=O(p^r)$ for some $r\le 1/4$; $\max_{i\le p}|\mu_i|/\sqrt{\sigma_{ii}}\ge \sqrt{2\beta \log(p)/n}$ with $\beta>1/\min_i(\sigma_{ii}\omega_{ii})+\varepsilon$ for some $\varepsilon>0$. Then, as $p\to\infty$
$$\P_{H_1}(\phi_\a(\Omega))\to 1,
$$
where $k_p$ is the number of non-zero entries of $\bm\mu$.
\begin{rmk}	
	Note that for  the statistic $T_{CLX}$, one can use any consistent estimator of $\Si^{-1}$ in the sense of the $L_1$-norm and infinity norm with at least a logarithmic rate of convergence. 
\end{rmk}

\subsection{MANOVA and Contrasts: more than two samples}

In this subsection, we consider the problem of testing the equality of several high-dimensional mean vectors, which is also called the multivariate analysis of variance (MANOVA) problem. This problem is to test the hypothesis
\begin{align}\lb{h}
	H_0:\bm\mu_1=\cdots=\bm\mu_k \quad \mbox{vs}\quad H_1:\exists i\neq j,~~ \bm\mu_i\neq\bm \mu_j.
\end{align} 
For samples that are drawn from a normal distribution family, the MANOVA  problem in a high-dimensional setting has been considered widely in the literature. For example, among others, Tonda and Fujikoshi \cite{TondaF04A} obtained the asymptotic null distribution of the likelihood ratio test; Fujikoshi \cite{Fujikoshi04M}  found the asymptotic null distributions for the Lawley-Hotelling trace and the Pillai trace statistics; and Fujikoshi et al \cite{FujikoshiH04A} considered the
Dempster trace test, which is based on the ratio of the trace of the between-class sample covariance matrix to the trace of the within-class sample covariance matrix.    Instead of investigating the ratio of the traces of the two sample matrices, Schott \cite{Schott07S} proposed a test statistic based on the difference between the traces.
Next, we introduce three NTM statistics that are the improvements on $T_{SD}$, $T_{CQ}$ and $T_{CLX}$. 

Recently, Srivastava and Kubokawa 
\cite{SrivastavaK13Ta} proposed a test statistic for testing the equality of the mean vectors of several groups with a common unknown non-singular covariance matrix.  Denote by ${\bf1}_r=(1,\dots,1)'$ an $r$-vector with all entries 1, and define
\begin{eqnarray*}
	\Y&=&(\X_{11},\dots,\X_{1n_1},\dots,\X_{k1},\dots,\X_{kn_k}),\\
\	\L&=&(\I_{k-1},-{\bf 1}_{k-1})_{(k-1)\times k}
\end{eqnarray*}
and $$\bfE=\left(
\begin{array}{ccc}
{\bf1}_{n_1} & \bf0 & \bf0 \\
\bf0 & {\bf1}_{n_2} & \bf0  \\
\vdots & \vdots & \vdots \\
\bf0 & \bf0 & {\bf1}_{n_k} \\
\end{array}
\right)_{n\times k}.
$$
Then, Srivastava and Kubokawa proposed the following test statistic:
\begin{align*}
	T_{SK}=\frac{{\rm tr}(\B\D_\S^{-1})-(n-k)p(k-1)(n-k-2)^{-1}}{\sqrt{2c_{p,n}(k-1)({\rm tr}\R^2-(n-k)^{-1}p^2)}},
\end{align*}
where $\B=\Y'\bfE(\bfE'\bfE)^{-1}\L'\bm [\L(\bfE'\bfE)^{-1}\L']^{-1}\L(\bfE'\bfE)^{-1}\bfE'\Y$, $\D_\S=Diag[(n-k)^{-1}\Y(\I_n-\bfE(\bfE'\bfE)^{-1}\bfE')\Y]$, $\R=\D_\S^{-1/2}\Y(\I_n-\bfE(\bfE'\bfE)^{-1}\bfE')\Y\D_\S^{-1/2}$ and $c_{p,n}=1+{\rm tr}(\R^2)/p^{3/2}.$ Note that $Diag[\A]$ denotes the diagonal matrix consisting of the diagonal elements of the matrix $\A$. Under the null hypothesis and the condition  $n\asymp p^\de$ with $\de>1/2$, $T_{SK}$ is asymptotically distributed as $N(0,1)$. Thus, as $n,p\to\infty$,
\begin{align*}
	P_{H_0}(T_{SK}>\xi_\a)\to\Phi(-\xi_\a).
\end{align*}
Hence, by comparing the results presented in \cite{SrivastavaD08T} and \cite{Srivastava09T}, it is easy to see that $c_{p,n}$ may be removable under certain conditions. 

Motivated by Chen and Qin \cite{ChenQ10T}, Hu et al. \cite{HuB14T} proposed a test for the MANOVA problem, that is, 
\begin{eqnarray*}
	T_{HB}&=&\sum_{i<j}^k(\bar {\X}_i-\bar{\bm  X}_j)'(\bar{\bm  X}_i-\bar {\X}_j)-(k-1)\sum_{i=1}^kn_i^{-1}{\rm tr}\S_i\\
	&=&(k-1)\sum_{i=1}^k\frac1{n_i(n_i-1)}\sum_{k_1\ne k_2}\X_{ik_1}'\X_{ik_2}-\sum_{i<j}^k\frac2{n_in_j}\sum_{k_1,k_2}\X_{ik_1}'\X_{jk_2}.
\end{eqnarray*}
When $k=2$, clearly, $T_{HB}$ reduces to Chen and Qin's test statistic. It is also shown that 
as $p\to\infty$ and $n\to\infty$,
\begin{align*}
	\frac{T_{HB}-\sum_{i<j}^k\| \mu_i- \mu_j\|^2}{\sqrt{Var(T_{HB})}}\toD N(0,1).
\end{align*}
To perform the test, a ratio-consistent estimator of $Var(T_{HB})$ for the MANOVA test is proposed in the paper. 

Cai  and Xia \cite{CaiX14H} also applied their idea to the MANOVA case under the homogeneous covariance assumption using the following test  statistic: 
\begin{align*}
	T_{CX}=\max_{1\leq i\leq p}\sum_{1\le j< l\leq k}\frac{n_in_j}{n_i+n_j}\left\{\frac{{\mathcal X}_{jli}^2}{\hat b_{ii}}\right\},
\end{align*}
where $\widehat{\Si^{-1}}(\bar X_j-\bar X_l):=({\mathcal X}_{jl1},\dots,{\mathcal X}_{jlp})'$; $\widehat{\Si^{-1}}:=(\o_{ij})_{p\times p}$ is a consistent estimator, e.g.  the constrained
$l_1$-minimization for the inverse matrix estimate of $\Si^{-1}$; and $\hat b_{ii}$ are the diagonal elements of the matrix $\hat{\B}$, which is defined by 
$$
\hat{\B}=\frac1{\sum n_i-k}\sum_{i=1}^k\sum_{j=1}^{n_i}\widehat\Si^{-1}(\X_{ij}-\bar{\X}_i)(\X_{ij}-\bar{\X}_i)'\widehat\Si^{-1}.
$$
To introduce the theory of Cai and Xia's test, let 
$${\mathcal Y}_i=\frac{1}{\hat\si_{ii}}\(\sqrt\frac{n_1n_2}{n_1+n_2}(\bar X_1-\bar X_2)_i,\dots, \sqrt\frac{n_{k-1}n_k}{n_{k-1}+n_k}(\bar X_{k-1}-\bar X_k)_i\)_{\frac{k(k-1)}{2}\times 1}',$$
where $\hat\si_{ii}$ is the estimate of the $(i,i)$-entry of the covariance matrix $\Si$. Let $\Si_{\cal Y}: \varrho\times \varrho$, $\varrho=\frac{k(k-1)}{2}$, be the covariance matrix of ${\mathcal Y}_i$. Let  $\lambda_{\cal Y}^2$ be the largest eigenvalue of $\Si_{\cal Y}$, and let $d$ be the dimension of the eigenspace of $\lambda_{\cal Y}^2$. Let  $\lambda_{\cal Y,i}^2: 1\le i\le\varrho$ be the eigenvalues 
of $\Si_{\cal Y}$ arranged in descending order. 

Under the null hypothesis $H_0$ and some regularity conditions on the population covariance matrix, for any $x\in R$, as $\min\{n,p\}\to \infty$,
\begin{align*}
	\P_{H_0}(T_{CX}-2\lambda_{\cal Y}^2\log(p)-(d-2)\lambda_{\cal Y}^2\log\log(p)\leq x)\to \exp\(-\Gamma^{-1}\(\frac{d}{2}\)H(\Si)\exp\(-\frac{x}{2\lambda_{\cal Y}^2}\)\),
\end{align*}
where $\Gamma$ is the gamma function, and $H=\prod_{i=d+1}^\varrho(1-\frac{\lambda_{\cal Y,i}^2}{\lambda_{\cal Y}^2})^{-1/2}$.
Similar to the two-sample location problem, they also established a theorem to evaluate the consistency of their test:

Suppose that  $C^{-1}_0\le \lambda_{\min}(\Si)\le \lambda_{\max}(\Si)\le C_0$ for some constant $C_0>1$. If 
$k_p=\max_{j<l\le k} \sum_{i=1}^pI(\mu_j-\mu_l\ne 0)=o(p^r)$, for some $r < 1/4$ and $\max_{i} \|\bm \delta_i\|_2/\sqrt{\si_{ii}}
\ge \sqrt{2\si^2\beta \log p}$ with some $\beta \ge  1/(\min_i \si_{ii}\omega_{ii}) + \varepsilon$ for some
constant $\varepsilon > 0$, then, as $p\to \infty$,
$$
\P_{H_1} (\phi_\a(\Omega)= 1) \to 1,
$$
where $\bm\delta_i=(\mu_{1i}-\mu_{2i},\cdots,\mu_{k-1,i}-\mu_{k,i})'$.

\subsection{Some related work on the tests of high-dimensional locations }
Chen et al. \cite{ChenP11R} proposed another statistic: 
\begin{align*}
	T_{RHT}=\bar{\X}'(\S+\la \I)^{-1}\bar{\X}, ~ \mbox{for $\la>0$,}
\end{align*}
which is called the regularized Hotelling $T^2$ test.  The idea is to employ the technique of ridge regression to stabilize the inverse of the sample covariance matrix given in (\ref{reht}).
Assuming that the underlying distribution is normally distributed, it is proven that under the null hypothesis, for any $\la >0$, as $p/n_1\to y\in(0,\infty)$ 
\begin{align*}
	\frac{\sqrt{p}\(nT_{RHT}/p-\frac{1-\la m(
			\la)}{1-p(1-\la m(\la))/n}\)}{\frac{1-\la m(\la)}{(1-p/n+p\la m(\la)/n)^3}-\la \frac{m(\la)-\la m'(\la)}{(1-p/n+p\la m(\la)/n)^4}}\toD N(0,1),
\end{align*} 
where $m(\la)=\frac{1}{p}{\rm tr}(\S+\la \I)^{-1}$ and $m'(\la)=\frac{1}{p}{\rm tr}(\S+\la \I)^{-2}$. They also give an asymptotic approximation method for selecting 
the tuning parameter $\la$ in the  regularization.  Recently, based on a supervised-learning strategy, Shen and Lin \cite{ShenL15A} proposed a statistic to select an optimal subset of features to maximize the asymptotic power of the Hotelling $T^2$ test.

The Random Projection was first proposed by Lopes et al. \cite{LopesJ11M} and was further discussed in later studies \cite{Thulin14H,WeiL13D,JacobN12D,ZhangP16H}.  For Gaussian data, the procedure projects high-dimensional data onto random subspaces of relatively low-dimensional spaces to allow the traditional Hotelling $T^2$ statistic to work well. 
This method can be viewed as a two-step procedure. First, a single random projection is drawn, and it is then used to map the samples from the high-dimensional space to a low-dimensional
space. Second, the Hotelling $T^2$ test is applied to a new hypothesis-testing problem in the projected space. A decision is then returned to the original problem by simply rejecting $H_0$ whenever the Hotelling test rejects it in the projected spaces.

Some other related work on tests of high-dimensional locations can be found in  \cite{HyodoN14O,BiswasG14N,TouloumisT15T,ChakrabortyC15W,MondalB15H,FengS,FengS15N}, which we do not discuss at length in this paper.

\section{NTM on covariance matrices}

\subsection{One-sample scatter test}

The standard test for scatters is to test the hypothesis $H_0:\Si=\Si_0$ vs $H_1: \Si\ne \Si_0$. Because $\Si_0$ is known, one can multiply $\Si_0^{-1/2}$ by the data set and then change the test to the simpler hypothesis $H_0:\Si=\I_p$. The classical test for this hypothesis is the well-known likelihood ratio, which can be found in any standard textbook, such as Anderson \cite{Anderson03I}. The likelihood ratio test statistic is given by 
$$
T_{LR}={\rm tr}{\S}-\log \det({\S})-p.
$$
When $p$ is fixed, the test based on $T_{LR}$ has many optimalities, such as unbiasedness, consistency, and being invariant under affine transformation. However, similar to the Hotelling $T^2$ test, it has a fatal defect in that it is not well defined when $p$ is larger than $n-1$. When $p$ is large but smaller than $n-1$, the null distribution is not simple to use, even under normality. The popularly used option is the Wilks theorem. However, when $p$ is large, the Wilks theorem introduces a very serious error to the test because its size tends to 1 as $p$ tends to infinity. A correction to the likelihood ratio test based on random matrix theory can be found in \cite{BaiJ09C}. However, when $p$ is large, especially when $p/n$ is close to 1, we believe that the asymptotic power will be low, much as occurs for the $T^2$ test. The idea of NTM can also be applied to this hypothesis. Now, we first introduce the work by Ledoit  and Wolf \cite{LedoitW02S}.

Ledoit and Wolf considered two hypotheses: $H_{01}:\Si=\I_p$ and $H_{02}:\Si=a\I_p$ with $a>0$ unknown. Based on the idea of the Nagao test (see Nagao (1973) \cite{Nagao73S}), they proposed two test statistics:
$$
V=\frac1p{\rm tr}(\S-\I_p)^2\ \ \mbox{ and }\ \ U=\frac1p{\rm tr}\(\frac{\S}{\frac1p{\rm tr}\S}-\I_p\)^2,
$$
which can be viewed from the perspective of NTM as considering $\S$ and  $\frac1p{\rm tr}\S$ to be the estimators of parameters $\Si$ and $a$ in the target functions
\begin{align}\lb{tag1}
h(\Si)=\frac1p{\rm tr}(\Si-\I_p)^2\ \ \mbox{ and }\ \ h(\Si, a)=\frac1p{\rm tr}\(\frac{\Si}{a}-\I_p\)^2
\end{align} , respectively. Note that under the null hypothesis, $a=\frac1p{\rm tr}\Si$.
They studied the asymptotic properties of $U$ and $V$ in the high-dimensional setting where $p/n\to c\in (0,\infty)$ and found that 
$U$, for the hypothesis of sphericity, is robust against large $p$, even larger than $n$. 
However, because $V$ is not consistent against every alternative, they proposed a new test statistic:
$$
W=\frac1p{\rm tr}(\S-\I_p)^2-\frac{p}{n}\left[\frac1p{\rm tr}\S\right]^2+\frac{p}{n}.
$$
Under normality and the assumptions that 
\begin{eqnarray*}
	\frac1p{\rm tr}(\Si)&=&\a\\
	\frac1p{\rm tr}(\Si-\a\I)^2&=&\delta^2\\
	\frac1p{\rm tr}(\bm \Si^j)&\to& \nu_j<\infty, \mbox{ for } j=3,4,
\end{eqnarray*} 
They proved the following: 
\begin{description}
	\item[(i)] The law of large numbers
	\begin{eqnarray*}
		\frac1p{\rm tr}(\S)&\toP& \alpha\\
		\frac1p{\rm tr}(\S^2)&\toP& (1+c)\a^2+\delta^2;
	\end{eqnarray*}
	\item[(ii)] The CLT, if $\delta=0$
	\begin{eqnarray*}
		&&n\begin{bmatrix}\frac1p{\rm tr}(\S)-\a\cr \frac1p{\rm tr}(\S)^2-\frac{n+p+1}{n}\a^2\cr \end{bmatrix}\\
		&\toD& N\left(\begin{bmatrix}0\cr 0\cr\end{bmatrix},\begin{bmatrix}\frac{2\a^2}{c}&4\left(1+\frac1c\right)\a^3\cr4\left(1+\frac1c\right)\a^3&4\left(\frac2c+5+2c\right)\a^4\cr\end{bmatrix}\right).
	\end{eqnarray*}
\end{description}
	Based on these results, they derived that $nU-p\toD N(1,4)$. 
The inconsistency of the test based on $V$ can be seen from the following facts.
When $p$ is fixed, by the law of large numbers, we have $\S\to \Si=\I_p$, and hence, $V=\frac1p {\rm tr }(\S-\I_p)^2\toP 0$. However, when $p/n\to c>0$, we have
\begin{eqnarray*}
	V&=&\frac1p{\rm tr}(\S)^2-\frac{2}p{\rm tr}(\S)+1\\
	&\toP& (1+c)\a^2+\delta^2-2\a+1=c\a^2+(\a-1)^2+\delta^2.
\end{eqnarray*}
Because the target function is $\frac1p{\rm tr}(\Si-\I_p)^2=(\a-1)^2+\delta^2$, the null hypothesis can be regarded as $(\a-1)^2+\delta^2=0$, and the alternative can be considered as $(\a-1)^2+\delta^2>0$. However, the limit of $V$ has one more term, $c\a^2$, which is positive. Therefore, the test $V$ is not consistent. 
In fact, it is easy to construct a counterexample based on this limit: set 
$$
c\a+(1-\a)^2+\delta^2=c.
$$
When $\delta=0$, the solution to the equation above is $\a=\frac{1-c}{1+c}$. Accordingly, the limit of $V$ is the same for the null $\a=1$ and the alternative $\a=\frac{1-c}{1+c}$. 

For $W$, we have 
$$
W\toP c\a^2+(\a-1)^2+\delta^2-c\a^2+c=c+(\a-1)^2+\delta^2.
$$
When $p$ is fixed, they proved that as $n\to\infty$,
$$
\frac{np}2W\toP \chi^2_{p(p+1)/2}
$$
or equivalently,
$$
nW-p\toP\frac2{p}\chi^2_{p(p+1)/2}-p
$$
When $p\to\infty$, the right-hand side of the above tends to $N(1,4)$, which is the same as the limit when $p/n\to c$. This behavior shows that the test based on $W$ is robust against $p$ increasing. Chen et al \cite{ChenZ10T} extended the work to the case without normality assumptions.

Now, the target functions \eqref{tag1} can be rewritten as 
\begin{align*}
h_1(\Si)=\frac1p{\rm tr}(\Si-\I_p)^2=\frac1p{\rm tr}\Si^2-\frac2p{\rm tr}\Si+1
\end{align*}
and 
\begin{align*}
h_2(a, \Si)=\frac1p{\rm tr}\(\frac{\Si}{a}-\I_p\)^2=\frac{\frac1p{\rm tr}\Si^2-(\frac1p{\rm tr}\Si)^2}{(\frac1p{\rm tr}\Si)^2}.
\end{align*}
 Then, under the normality assumption, Srivastava \cite{Srivastava05S} gave the unbiased and consistent estimators of these parameters in the previous target functions, which are as follows:
 \begin{align*}
\widehat{\frac1p{\rm tr}\Si}=\frac{1}{p}{\rm tr}\S\ \ \mbox{and}\ \ \widehat{\frac1p{\rm tr}\Si^2}=\frac{(n-1)^2}{p(n-2)(n+1)}\({\rm tr}\S^2-\frac{1}{n-1}({\rm tr}\S)^2\).
 \end{align*}
 Based on these estimators, he proposed the test statistics 
 \begin{align*}
 T_{S1}=\widehat{\frac1p{\rm tr}\Si^2}-2\widehat{\frac1p{\rm tr}\Si}+1 \ \ \mbox{and}\ \ T_{S2}=\frac{\widehat{\frac1p{\rm tr}\Si^2}-(\widehat{\frac1p{\rm tr}\Si})^2}{(\widehat{\frac1p{\rm tr}\Si})^2},
 \end{align*} 
and  proved that under the assumption $n\asymp p^\de, 0<\de
 \leq 1$, as $\{n,p\}\to\infty$, we have asymptotically,
 \begin{align*}
 \frac{n}{2}\(T_{S1}-\frac1p{\rm tr}(\Si-\I_p)^2\)\sim N(0,\tau_1^2)
 \end{align*}
 and 
  \begin{align*}
  \frac{n}{2}\(T_{S2}-\frac{\frac1p{\rm tr}\Si^2-(\frac1p{\rm tr}\Si)^2}{(\frac1p{\rm tr}\Si)^2}\)\sim N(0,\tau_2^2),
  \end{align*}
 where $\tau_1^2=\frac{2n}{p}(\a_2 -2\a_3 +\a_4)+\a_2^2$, $\tau_2^2=\frac{2n(\a_4\a^2_1 - 2\a_1\a_2\a_3 + a^3_2)}{p\a_1^6}+\frac{\a_2^2}{\a_1^4}$ and $\a_i=\frac1p{\rm tr}\Si^i$. Thus, under the null hypothesis, one can easily obtain \begin{align*}
  \frac{n}{2}T_{S1}\toD N(0,1)\ \ \mbox{and }\ \   \frac{n}{2}T_{S2}\toD N(0,1).
 \end{align*}
 
   Later, Srivastava and Yanagihara \cite{SrivastavaY10T}  and  Srivastava et al. \cite{SrivastavaY14T} extended this work to the cases of two or more population covariance matrices and without normality assumptions, respectively. 
Furthermore,  Cai and Ma \cite{CaiM13O} showed that $\T_{S1} $ is rate-optimal over this asymptotic regime, and Zhang et al. \cite{ZhangP13T} proposed the empirical likelihood ratio test for this problem.

\subsection{ Li and Chen's test based on unbiased estimation of target function}

Li and Chen (2012) \cite{LiC12T} considered the two-sample scatter problem, that is, testing the hypothesis $H_0: \Si_1=\Si_2$. They choose the target function as $h(\Si_1,\Si_2)={\rm tr}(\Si_1-\Si_2)^2$. 
They selected the test statistic by the unbiased estimator of $h(\Si_1,\Si_2)$ as 
$$
T_{LC}=A_{n_1}+A_{n_2}-2C_{n_1n_2}
$$
where 
\begin{eqnarray*}
	A_{n_h}&=&\frac1{(n_h)_2}\sum_{i\ne j}(\X_{hi}'\X_{hj})^2-\frac{2}{(n_h)_3}\sum_{i,j,k\atop distinct}
	\X_{hi}'\X_{hj}\X_{hj}'\X_{hk}\\
	&&+\frac1{(n_h)_4}\sum_{i,j,k,l\atop distinct}\X_{hi}'\X_{hj}\X_{hk}'\X_{hl},\\\\
	C_{n_1n_2}&=&\frac1{n_1n_2}\sum_{i, j}(\X_{1i}'\X_{2j})^2-\frac{1}{n_2(n_1)_2}\sum_{i\neq k}\sum_{j}
	\X_{1i}'\X_{2j}\X_{2j}'\X_{1k}\\
	&&-\frac{1}{n_1(n_2)_2}\sum_{i\neq k}\sum_{j}
	\X_{2i}'\X_{1j}\X_{1j}'\X_{2k}
	+\frac1{(n_1)_2(n_2)_2}\sum_{i\neq j}\sum_{k\neq l}\X_{1i}'\X_{2j}\X_{1k}'\X_{2l}.
\end{eqnarray*}
Under the conditions A1 and A2 and for any $i,j,k,l\in\{1,2\}$,
$$
{\rm tr}( \Si_i\Si_j\Si_k\Si_l)=o({\rm tr}(\Si_i\Si_j){\rm tr}(\Si_k\Si_l)),
$$
we have
\begin{eqnarray*}
	Var (T_{LC})&=&\sum_{i=1}^2\bigg[\frac4{n_i^2}{\rm tr}^2\Si_i^2+\frac{8}{n_i}{\rm tr}(\Si_i^2-\Si_1\Si_2)^2\\
	&& +\frac{4}{n_i}{\rm tr}(\Gamma_i'(\Si_1-\Si_2)\Gamma_i\circ \Gamma_i'(\Si_1-\Si_2)\Gamma_i)\bigg]\\
	&&+\frac{8}{n_1n_2}{\rm tr}^2(\Si_1\Si_2),
\end{eqnarray*}
where $\A\circ \B=(a_{ij}b_{ij})$ denotes the Hadamard product of matrices $\A$ and $\B$.

Li and Chen \cite{LiC12T} proved that
$$
\frac{T_{CL}-{\rm tr}(\Si_1-\Si_2)^2}{\sqrt{Var(T_{LC})}}\toD N(0,1).
$$
Li and Chen selected $\widehat{\sqrt{Var(T_{CL})}}:= \frac{2}{n_1}A_{n_1}+\frac{2}{n_2}A_{n_2}$, which is a ratio-consistent estimator of $\sqrt{Var(T_{LC})}$ under $H_0$. Therefore, the test rejects $H_0$ if 
$$
T_{LC}>\xi_\a\(\frac{2}{n_1}A_{n_1}+\frac{2}{n_2}A_{n_2}\).
$$
\begin{rmk}
	In \cite{LiC12T},  Li and Chen also considered the test for the covariance between two sub-vectors, i.e., testing the hypothesis $H_0: ~\Sigma_{1,12}=\Sigma_{2,12}$, where  $\Sigma_{i,12}$ is the off-diagonal blocks of $\Sigma_i$. As the test statistic is similar, we omit the details here.
\end{rmk}
\subsection{ Cai et al's maximum difference test}

Cai et al \cite{CaiL13T} also applied their maximum elements of the difference of two sample covariance matrices to test the hypothesis of the equality of the two population covariances. They defined their test statistic as follows:
$$
M_n=\max_{1\le i\le j\le p}M_{ij}=\max_{1\le i\le j\le p}\frac{(s_{ij1}-s_{ij2})^2}{\hat\theta_{ij1}/n_1+\hat\theta_{ij2}/n_2},
$$
where $s_{ijl}$ is the $(i,j)$-th element of the sample covariance of the $l$-th sample, and 
$$
\hat\theta_{ijl}=\frac1{n_l}\sum_{k=1}^{n_l}\left[(X_{kil}-\bar{X}_{il})(X_{kjl}-\bar{X}_{jl})-s_{ijl}\right]^2
$$
$1\le i\le j\le p$ and $l=1,2$. Here, $\hat\theta_{ijl}$ can be considered an estimator of the variance of 
$s_{ijl}$. Then, they defined the test by
$$
\phi_\a=I(M_n>q_\a+4\log p-\log\log p).
$$
where $q_\a$ is the upper $\a$ quantile of the Type I extreme value distribution with the c.d.f. 
$$\exp\Big(-\frac1{\sqrt{8\pi}}\exp(-\frac{x}2)\Big),$$
and therefore
$$
q_\a=-\log(8\pi)-2\log\log(1-\a)^{-1}.
$$
Under sparse conditions on the difference of the population covariances $\Si_1-\Si_2$ and certain distributional conditions, they proved that for any $t\in R$
$$
\P(M_1-4\log p+\log\log p\le t)\to \exp\left(-\frac1{\sqrt{8\pi}}\exp\left(-\frac{t}{2}\right)\right).
$$
As expected, Cai et al's test is powerful when the difference of the two population covariances is sparse, and it thus compensates somewhat for Li and Chen's test.

\subsection{Testing the structure of the covariance matrix  }
In this subsection, we will consider another important test problem, namely, testing the structure of the covariance matrix. First, we review the test hypothesis that the covariance matrix $\Sigma$ is banded. That is, the variables have nonzero correlations only up to a certain lag $\tau\geq1$. To elaborate, we denote $\Sigma=(\sigma_{ij})_{p\times p}$  and consider the following test hypotheses: \begin{align}\lb{teststr1}
H_0:~\sigma_{ij}=0,~~\mbox{for all}~ |i-j|\geq \tau~~\mbox{v.s.} ~~ H_1:~\sigma_{ij}\neq 0,~~\mbox{for some}~ |i-j|\geq \tau,
\end{align}
or, equivalently, 
\begin{align*}
H_0: ~\Sigma=B_\tau(\Sigma)~~\mbox{v.s.} ~~ H_1:~\Sigma\neq B_\tau(\Sigma),
\end{align*}
where $B_\tau(\Sigma)=(\sigma_{ij}I(|i-j|\leq \tau)).$
From the perspective of NTMs, one can also choose the target functions by the Euclidean distance and the Kolmogorov distance, which are the main concepts of the tests proposed by Qiu and Chen \cite{QiuC12T} and 
Cai and Jiang \cite{CaiJ11L}, respectively.

For $\tau+1\leq q\leq p-1$ and ${ \bm \mu}={\bf 0}$, let $$\widehat{\sigma_{ll+q}^2}=\frac{1}{n(n-1)}\sum_{i\neq j}X_{li}X_{(l+q)i}X_{lj}X_{(l+q)j}-\frac{2}{n(n-1)(n-2)}\sum_{i,j,k\atop distinct}X_{li}X_{(l+q)j}X_{lk}X_{(l+q)k}$$$$+\frac{1}{n(n-1)(n-2)(n-3)}\sum_{i,j,k,m\atop distinct}X_{li}X_{(l+q)j}X_{lk}X_{(l+q)m}.$$
By denoting $T_{QC}^{\tau}:=2\sum_{q=k+1}^{p-1}\sum_{l=1}^{p-q}\widehat{\sigma_{ll+q}^2}$, one can easily check that $T_{QC}^{\tau}$ is an unbiased estimator of ${\rm tr}(\Sigma-B_\tau(\Sigma))$. Under the 
assumptions $\tau=o(p^{1/4})$, (A1), (A2) and certain conditions on the eigenvalues of $\Sigma$, Qiu and Chen \cite{QiuC12T} showed that 
under the null hypothesis,
\begin{align*}
\frac{nT_{QC}^\tau}{V_{n\tau}}\toD N(0,4),
\end{align*}
and the power function asymptotically satisfies 
\begin{align*}
\beta_{QC}=\P({nT_{QC}^\tau}/{V_{n\tau}}\geq 2\xi_\a|\Sigma\neq B_k(\Sigma))
\simeq \Phi\left(\frac{2\xi_\a V_{n\tau}}{nv_{n\tau}}-\delta_{n\tau}\right) \geq \Phi\left(\frac{\xi_\a V_{n\tau}}{{\rm tr}(\Sigma^2)}-\delta_{n\tau}\right),
\end{align*}
where $V_{n\tau}=\sum_{l=1}^p\widehat {\sigma_{ll}^2}+2\sum_{q=1}^{\tau}\sum_{l=1}^{p-q}\widehat{\sigma_{ll+q}^2}$, $v_{n\tau}^2=4n^{-2}{\rm tr}^2(\Sigma^2)+8n^{-1}{\rm tr}(\Sigma(\Sigma-B_\tau(\Sigma)))^2+4n^{-1}\Delta
{\rm tr}[(\Gamma'(\Sigma-B_\tau(\Sigma))\Gamma)\circ (\Gamma'(\Sigma-B_\tau(\Sigma))\Gamma)]$ and $\delta_{n\tau}={\rm tr}(\Sigma-B_\tau(\Sigma)^2/v_{n\tau}$.

We can also rewrite the test hypothesis \eqref{teststr1} as 
$$H_0:~\rho_{ij}=0,~~\mbox{for all}~ |i-j|\geq \tau~~\mbox{v.s.} ~~ H_1:~\rho_{ij}\neq 0,~~\mbox{for some}~ |i-j|\geq \tau,$$ 
where $\rho_{ij}$ is the population correlation coefficient between two random variables $X_{1i}$ and $X_{1j}$. Cai and Jiang \cite{CaiJ11L} proposed a test procedure based on the largest magnitude
of the off-diagonal entries of the sample correlation matrix 
$$T_{CL}^{\tau}=\max_{|i-j|\geq\tau}|\hat{\rho
	}_{ij}|,$$ where $
\hat{\rho
}_{ij}$ is the sample correlation coefficient.
They showed that under the assumptions $\log p=o(n^{1/3})\to\infty$ and $\tau=o(p^\epsilon)$ with $\epsilon>0$, for any $t\in R$,
\begin{align*}
	\P(n(T_{CL}^{\tau})^2-4\log p+\log\log p\leq t)\to \exp\left(-\frac1{\sqrt{8\pi}}\exp\left(-\frac{t}{2}\right)\right).
\end{align*}
By implication, one can reject the null hypothesis whenever \begin{align*}
(T_{CL}^{\tau})^2\geq n^{-1}[4\log p-\log\log p-\log(8\pi)-2\log\log(1-\a)^{-1}]
\end{align*}
with asymptotical size $\a$.

\begin{rmk}
	If $\tau=1$ and under the normal assumption, then the test hypothesis is (\ref{teststr1}), also known as testing for complete independence, which was first considered by Schott in 2005 \cite{Schott05T} for a high-dimensional random vector and using the Euclidean distance of the sample correlation matrix. 
	\end{rmk}
\begin{rmk}
	Peng et al. \cite{PengC16M} improved  the power of the test $T_{QC}^{\tau}$ by employing the banding estimator for the covariance matrices.
	Zhang et al. \cite{ZhangP13T} also gave the empirical likelihood ratio test procedure for testing whether the population covariance matrix has a banded structure.
\end{rmk}

\section{Conclusions and Comments}

All of the NTM procedures show that most classical procedures in multivariate analysis are less powerful in some parameter settings when the dimension of data is large. Thus, it is necessary to develop new procedures to improve the classical ones. However, all of the NTM procedures developed to date require additional conditions on the unknown parameters to guarantee the optimality of the new procedures; e.g., all procedures based on asymptotically normal estimations require that the eigenstructure of population covariance matrix should not be too odd, and all NTMs based on the Kolmogorov distance require the sparseness of the known parameters. Therefore, there is a strong need to develop data-driven procedures that are optimal in most cases.

\Acknowledgements{ The authors would like to thank the referees for their constructive comments, which led to a substantial improvement of the paper. J. Hu was partially supported by  the National Natural Science Foundation of China (Grant No.  11301063),  Science and Technology Development foundation of Jilin (Grant No.  20160520174JH),  and Science and Technology Foundation of  Jilin during the ``13th Five-Year Plan"; and Z. D. Bai was partially supported by the National Natural Science Foundation of China (Grant No.  11571067). }


\begin{thebibliography}{10}
\bibitem{Fujikoshi04M}
	Y.~Fujikoshi.
	\newblock {Multivariate analysis for the case when the dimension is large
		compared to the sample size}.
	\newblock {\em Journal of the Korean Statistical Society}, 33(1):1--24, 2004.
\end{thebibliography}
%

\end{document}